\documentclass[11pt]{article}
\usepackage{amssymb, amsmath,  amsfonts,dsfont, mathrsfs,euscript,eufrak,colonequals}
\usepackage{textcomp}
\usepackage{latexsym}
\usepackage{indentfirst}
\usepackage[normalem]{ulem}

\newtheorem{Theorem}{Theorem}
\newtheorem{Lemma}{Lemma}
\newtheorem{Proposition}{Proposition}
\newtheorem{Corollary}{Corollary}

\newcommand{\defeq}{\colonequals}

\newcommand{\dd}{\,d}
\newcommand{\id}{\mathds{1}}

\newcommand{\Complex}{\mathbb C}
\newcommand{\R}{\mathbb R}
\newcommand{\Z}{\mathbb Z}
\newcommand{\N}{\mathbb N}

\newcommand{\Realpart}{\mathrm{Re}\,}
\newcommand{\Imagpart}{\mathrm{Im}\,}

\newcommand{\BV}{\boldsymbol{V}}

\newcommand{\e}{\varepsilon}

\newcommand{\cconv}{\Rightarrow}

\newcommand{\Ln}{\mathrm{Ln}}

\newcommand{\aeset}{\mathcal{T}}

\begin{document}
\author{A. A. Khartov\footnote{Laboratory for Approximation Problems of Probability,  Smolensk State University, 4 Przhevalsky st., 214000 Smolensk, Russia, e-mail: \texttt{alexeykhartov@gmail.com} }}
\title{On weak convergence of quasi-infinitely divisible laws}

\maketitle

\begin{abstract}
	We study a new class of so-called quasi-infinitely divisible laws, which is a wide natural extension of the well known class of  infinitely divisible laws through the L\'evy--Khinchine type representations. We are interested in criteria of weak convergence within this class. Under rather natural assumptions, we state assertions, which  connect a weak convergence of quasi-infinitely divisible distribution functions with one special type of convergence of their L\'evy--Khinchine spectral functions. The latter convergence is not equivalent to the weak convergence. So we complement known results by Lindner, Pan, and Sato (2018) in this field. 
\end{abstract}


\section{Introduction}

This paper is devoted to the questions concerning weak convergence within a new class of so-called quasi-infinitely divisible probability laws.

Let $F$ be a  distribution function of a probability law on the real line $\R$ with the characteristic function 
\begin{eqnarray*}
	f(t)\defeq\int_{\R} e^{itx} \dd F(x),\quad t\in\R.
\end{eqnarray*}
Recall that $F$ (and the corresponding law) is called \textit{infinitely divisible} if for every positive integer $n$ there exists a  distribution function $F_{1/n}$ such that $F=(F_{1/n})^{*n}$, where ``$*$'' is the convolution, i.e. $F$ is $n$-fold convolution power of $F_{1/n}$.  It is well known that $F$ is infinitely divisible if and only if the characteristic function $f$ is represented by \textit{the L\'evy--Khinchine formula}:
\begin{eqnarray}\label{def_ReprLevyKhinchine}
	f(t)=\exp\biggl\{it \gamma+\int_{\R} \bigl(e^{itx} -1 -\tfrac{it}{\tau} \sin(\tau x)\bigr)\tfrac{1+x^2 }{x^2} \dd G(x)\biggr\},\quad t\in\R,
\end{eqnarray}
with some $\tau>0$, \textit{shift parameter} $\gamma\in\R$, and with a  bounded non-decreasing \textit{spectral function} $G: \R \to \R$  that is assumed to be right-continuous with condition $G(-\infty)=0$ (throughout the paper, $G(\pm\infty)$ denote the limits at $\pm \infty$ correspondingly). We use $u\mapsto \tfrac{1}{\tau}\sin (\tau u)$ as the ``centering function'' in the integral in \eqref{def_ReprLevyKhinchine} following to Zolotarev \cite{Zolot1} and \cite{Zolot2}. If formula  \eqref{def_ReprLevyKhinchine} holds for some $\tau=\tau_0>0$, then it holds for any $\tau>0$, where $\gamma$ will depend on $\tau$, but $G$ will not. It is well known that \textit{the spectral pair} $(\gamma, G)$ is uniquely determined by $f$ and hence by $F$. The L\'evy--Khinchine formula plays a fundamental role in  probability theory, and it also has a lot of applications in related fields (see \cite{Appl} and \cite{Sato}).

It turns out that there exists a rather wide class of probability laws that are very similar to infinitely divisible laws. This class of so-called \textit{quasi-infinitely divisible laws} was introduced by Lindner and Sato \cite{LindSato}.  Following them, a distribution function $F$ (and the corresponding law) is called \textit{quasi-infinitely divisible} if its characteristic function $f$ admits the representation \eqref{def_ReprLevyKhinchine} with some \textit{shift parameter}  $\gamma\in\R$, \textit{spectral function} $G: \R \to \R$ of bounded variation on $\R$ (not necessarily monotone), and for some (any) $\tau>0$.  Here $G$ is assumed to be  right-continuous with condition $G(-\infty)=0$ as before and so  $f$ (and $F$) uniquely determines \textit{the spectral pair} $(\gamma, G)$ (see \cite{GnedKolm} p. 80). Observe that, due to the Jordan decomposition, we can represent $G(x)=G^+(x)-G^-(x)$, $x\in\R$, with some bounded non-decreasing functions $G^+$ and $G^-$ on $\R$.   Also we  can always write $\gamma=\gamma^+-\gamma^-$ with some numbers $\gamma^+$ and $\gamma^-$ from $\R$. Then it is clear that $f(t)=f^{+}(t)/f^{-}(t)$, $t\in\R$,  where $f^+$ and $f^-$ are characteristic functions of some two infinitely divisible distribution functions $F^+$ and $F^-$ with the spectral pairs $(\gamma^+,G^+)$ and $(\gamma^-,G^-)$ correspondingly, and so $F*F^-=F^+$. Starting from this point of view, it is rather natural to call distribution function $F$ (and the corresponding law) \textit{rationally infinitely divisible}. So every infinitely divisible law is quasi-infinitely divisible, but the converse is not true. There are a lot of interesting examples of quasi-infinitely divisible laws, which are not infinitely divisible (see \cite{GnedKolm} p. 82--83, \cite{LinOstr} p. 165, \cite{Lukacs} p. 123--124). Moreover, it seems that the class of quasi-infinitely divisible laws is essentially wider  than the class of infinitely divisible ones. In particular, it is clearly seen within discrete probability laws (see  \cite{AlexeevKhartov},  \cite{Khartov2}, \cite{KhartAlexeev}, and \cite{LindPanSato}).

Various forms of definition and the  first detailed analysis of the class of quasi-infinitely divisible laws on $\R$ was performed in \cite{LindPanSato}, the multivariate case is considered in the recent papers \cite{BergKutLind}, \cite{BergLind}, and \cite{Kutlu}. There are some results for discrete probability laws in this field (see \cite{AlexeevKhartov}, \cite{AlexeevKhartov2}, \cite{Khartov}, \cite{Khartov2}, and \cite{KhartAlexeev}) and for mixed laws (see \cite{Berger} and  \cite{BergerKutlu}).  It should be noted that quasi-infinitely divisible laws now have interesting applications in theory of stochastic processes (see  \cite{LindSato} and \cite{Pass}), number theory (see \cite{Nakamura}),  physics (see \cite{ChhDemniMou} and \cite{Demni}), and insurance mathematics (see  \cite{ZhangLiuLi}). 

We now focus on a weak convergence of quasi-infinitely divisible laws. Recall that, by definition, the sequence $(F_n)_{n\in\N}$ (where $\N$ is the set of positive integers) of distribution functions \textit{weakly converges} to a distribution function $F$ (we will write $F_n \xrightarrow{w} F$, $n\to\infty$) if 
\begin{eqnarray}\label{def_weakconv}
	\int_{\R} h(x)\dd F_n(x)\to \int_{\R} h(x)\dd F(x),\quad n\to\infty,
\end{eqnarray}
for any bounded continuous function $h:\R\to\R$. It is well known fact that this is equivalent to the following convergence:
\begin{eqnarray}\label{def_pointwiseconv}
	F_n(x)\to F(x), \quad n\to\infty, \quad \text{for any}\quad x\in S, 
\end{eqnarray}
where $S$ is an arbitrary dense subset of $\R$ and, in particular, it can be choosen as the set of all continuity points of $F$. The latter convergence is usually called \textit{weak} too (see \cite{Lukacs}). 

The weak convergence is also introduced for the class of real-valued functions of bounded variation on the real line (or for corresponding signed measures). Following Bogachev \cite{Bogachev}, it is analogously defined by formula \eqref{def_weakconv}, but instead of $F_n$, $n\in\N$, and $F$ we write some functions of bounded variation $G_n$, $n\in\N$, and $G$ correspondingly.  We will save the notation $G_n\xrightarrow{w} G$, $n\to\infty$, in this case. It should be noted that here weak convergence is not equivalent to the analog of convergence \eqref{def_pointwiseconv} with functions of bounded variation (see \cite{Bogachev} Section 1.4).   

There are rather general results by Lindner, Pan, and Sato in \cite{LindPanSato} concerning the weak convergence of quasi-infinitely divisible distribution functions. The authors state the conditions under which the weak convergence of distribution functions implies the weak convergence of the corresponding spectral functions together with the convergence of the shift parameters and vice versa. Namely, let $(F_n)_{n\in\N}$ be a sequence of quasi-infinitely divisible distribution functions and let $(\gamma_n, G_n)$ be the spectral pair of $F_n$ for every $n\in\N$. Let $F$ be a quasi-infinitely divisible distribution function with the spectral pair $(\gamma, G)$. Then the results from \cite{LindPanSato}  are in fact the following: 1) if $\gamma_n\to \gamma$ and $G_n\xrightarrow{w} G$, $n\to\infty$, then  $F_n\xrightarrow{w} F$, $n\to\infty$; 2) if we suppose $F_n\xrightarrow{w} F$, $n\to\infty$, then, under some assumptions on tightness and uniform boundedness for $(G_n)_{n\in\N}$, we have $\gamma_n\to \gamma$ and $G_n\xrightarrow{w} G$, $n\to\infty$. Here we omitted some details, the full formulation will be given in Section 3.

In this work we complement the  results by Lindner, Pan, and Sato. We  connect the weak convergence of quasi-infinitely divisible distribution functions with one type of convergence of their spectral functions. The latter convergence is a special modification of the convergence \eqref{def_pointwiseconv} (see the next section for details), and we think that it is more natural and explicit than the weak convergence for the class of functions of bounded variation.  A very similar convergence appeared in \cite{Bogachev} in Theorem 1.4.7., but we are not aware of existence of a definition for such convergence. So we introduce the necessary definition in Section 2.  We also show that the introduced convergence for function of bounded variations follows from the pointwise convergence of their Fourier--Stieltjes transforms under some natural assumptions.  This and other close propositions are key tools for  our main results devoted the weak convergence of quasi-infinitely divisible distribution functions. The main results are presented in Section 3. All proofs are provided in Section 4.

\section{Preliminaries and tools}

Let us consider the class of all functions $G: \R\to\R$ of bounded variations on $\R$.  Since we will interesting in the functions $G$, which generate the the Lebesgue--Stieltjes signed measures $\mu_G$,  we will focus only on  right-continuous functions $G$.  So for  the measures we will have  $\mu_G((a,b])=G(b)-G(a)$ for all $a,b\in\R$ and $a\leqslant b$. Recall that all intervals $(a,b]$ consist the generating semiring for $\mu_G$. So if there are two right-continuous functions $G_1$ and $G_2$ of bounded variations such that $G_2(x)=G_1(x)+C$, $x\in\R$, where $C\in\R$ is a constant, then the corresponding measures are  the same. Therefore we will consider only functions $G$ that satisfy $G(-\infty)=0$. 

Let $\BV$ denote the class of all functions  $G: \R\to\R$ of bounded variation on $\R$, which are right-continuous at every point $x\in\R$ and satisfy $G(-\infty)=0$. For every $G\in\BV$ its total variation on $\R$ will be denoted by $\|G\|$ and the total variation on $(-\infty,x]$ --- by $|G|(x)$, $x\in\R$. So we have 
\begin{eqnarray}\label{ineq_fbv_var}
	|G(x)|\leqslant |G|(x)\leqslant\|G\|,\quad\text{for any} \quad x\in\R,
\end{eqnarray}
and $ |G|(+\infty)=\|G\|$.

We now introduce a special type of convergence on the class $\BV$. Suppose that a whole sequence $(G_n)_{n\in\N}$ and a function $G$ are from $\BV$. We say that $(G_n)_{n\in\N}$ \textit{converges basically} to $G$, and write $G_n \cconv G$, $n\to\infty$, if    
 every its subsequence contain a further subsequence $(G_{n_k})_{k\in\N}$ such that
\begin{eqnarray*}
	G_{n_k}(x_2)-G_{n_k}(x_1)\to G(x_2)-G(x_1),\quad k\to\infty,
\end{eqnarray*}
for any $x_1, x_2\in\R$ except  at most a countable set, which in general depends on the choice of the subsequences.

Let us show that the basic convergence is equivalent to the weak convergence for distribution functions. Let $(F_n)_{n\in\N}$ be a sequence of distribution functions and let $F$ be a distribution function. Suppose that $F_n\xrightarrow{w}F$, $n\to\infty$. Then we have \eqref{def_pointwiseconv}, where $S$ is the set of all continuity points of $F$. Hence
\begin{eqnarray*}
	F_n(x_2)-F_n(x_1)\to F(x_2)-F(x_1),\quad n\to\infty, \quad \text{for all}\quad x_1,x_2 \in S.
\end{eqnarray*}
Since $\R\!\setminus\! S$ is at most countable set, we conclude that $F_n\cconv F$, $n\to\infty$, by definition. We now suppose that $F_n\cconv F$, $n\to\infty$. Let $(F_{n_k})_{k\in\N}$ be arbitrary subsequence of $(F_n)_{n\in\N}$ such that   
\begin{eqnarray}\label{conc_Fnk}
	F_{n_k}(x_2)-F_{n_k}(x_1)\to F(x_2)-F(x_1),\quad k\to\infty,
\end{eqnarray}
for any $x_1, x_2\in\R$ except  at most a countable set $D$. Let us fix $\e>0$ and choose $r_\e>0$ such that $\pm r_\e\in \R\!\setminus\! D$ and  $1-F(r_\e)+F(-r_\e)<\e$. We define $T_k(r)\defeq 1-F_{n_k}(r)+F_{n_k}(-r)$, $k\in\N$, $r>0$. Due to \eqref{conc_Fnk},  there exists $k_\e\in\N$ such that  $T_k(r_\e)<\e$ for all $k \geqslant k_\e$. Taking $r_\e$ greater to provide  $T_k(r_\e)<\e$ for all $k < k_\e$, we obtain $\sup_{k\in\N} T_k(r_\e)<\e$ because, due to monotonicity of every $F_{n_k}$, the inequality $T_k(r_\e)<\e$ still holds for all $k \geqslant k_\e$. Thus $\sup_{k\in\N} T_k(r)\to 0$, $r\to\infty$, and, in particular,  $\sup_{k\in\N} F_{n_k}(-r)\to0$, $r\to\infty$. Due to the latter, it is easy to check that \eqref{conc_Fnk} yields the convergence $F_{n_k}(x)\to F(x)$ for any $x\in\R$ except at most countable set $D$. Since  $\R\!\setminus\! D$ is dense subset of $\R$,  we have $F_{n_k}\xrightarrow{w} F$, $k\to\infty$.  Thus, according to  definition of a basic convergence,  every  subsequence of $(F_n)_{n\in\N}$ contains a further subsequence $(F_{n_k})_{k\in\N}$ that satisfies \eqref{conc_Fnk} and hence weakly converges to $F$. By the well known fact (see  \cite{Bill} p. 337), it means that the whole  sequence $(F_n)_{n\in\N}$ weakly converges to $F$.

The proved assertion can be generalized for bounded non-decreasing functions $F\in\BV$ and $F_n\in\BV$, $n\in\N$, but here the basic convergence $F_n\cconv F$, $n\to\infty$, must be taken  together with an additional condition that $F_n(+\infty)\to F(+\infty)$, $n\to\infty$ (see \cite{GnedKolm}, p. 39). It should be noted that basic and weak convergences  are not equivalent in a general case for functions from $\BV$. Indeed, the weak convergence implies the basic one that will follow from Theorem \ref{th_FStoFBV} below,  and also it is seen from Theorem~1.4.7. in \cite{Bogachev}. However, the converse is not true. The latter assertion is concluded from  the following simple examples.

\textbf{Example 1.}
Let us define $G_n(x)\defeq  \id_{n}(x)-\id_{n+1}(x)$, $x\in\R$, $n\in\N$. It is easily seen that $G_n\in \BV$, $n\in\N$, and $G_n(x)\to 0$, $n\to\infty$, for all $x\in\R$. Setting $G(x)\defeq 0$, $x\in\R$, we have the basic convergence $G_n\cconv G$, $n\to\infty$. Here $G_n(+\infty)=G(+\infty)=0$, and $\|G_n\|=2$, $n\in\N$. However, for the continuous and bounded function $x\mapsto \cos (\pi x)$, $x\in\R$, we conclude that $\int_{\R} \cos (\pi x) \dd G_n(x) \nrightarrow\int_{\R} \cos (\pi x) \dd G(x)=0$, $n\to\infty$. Indeed, 
\begin{eqnarray*}
	\int_{\R} \cos (\pi x) \dd G_n(x)=\cos (\pi n) - \cos (\pi (n+1))=(-1)^n-(-1)^{n+1}=2\cdot (-1)^n\nrightarrow 0,\quad n\to\infty.
\end{eqnarray*}
Thus $(G_n)_{n\in\N}$ doesn't weakly converge to $G$.\quad $\Box$\\

\textbf{Example 2.}
Let  $G_n(x)\defeq n \id_0(x)-n \id_{1/n^2}(x)$, $x\in\R$, $n\in\N$. So $G_n\in \BV$, $n\in\N$, and $G_n(x)\to 0$, $n\to\infty$, for all $x\ne 0$. We set $G(x)\defeq 0$, $x\in\R$, and we obtain that $G_n\cconv G$, $n\to\infty$. Observe that $\|G_n\|=2n\to\infty$, $n\to\infty$. Hence $(G_n)_{n\in\N}$ can not be weak convergent sequence, because, under the weak convergence, total variations  must be uniformly bounded (see Proposition 1.4.4. in \cite{Bogachev}, p. 22). Moreover, it even fails to hold that 
\begin{eqnarray}\label{def_vague}
	\int_{\R} h(x)\dd G_n(x)\to \int_{\R} h(x)\dd G(x),\quad n\to\infty,
\end{eqnarray}
for any continuous function $h$ with compact support. Indeed, let $h(x)\defeq\sqrt{x}$ for $x\in[0,1]$, $h(x)\defeq\sqrt{2-x}$ for $x\in[1,2]$, and $h(x)\defeq0$ for $x\notin [0,2]$. Obviously, the function $h$ satisfies the required properties.  So we have
\begin{eqnarray*}
	\int_{\R} h(x)\dd G_n(x)= h(0)\cdot n - h(1/n^2)\cdot n=0-\sqrt{1/n^2}\cdot n=-1,\quad \text{for every}\quad n\in\N,
\end{eqnarray*}
but $\int_{\R} h(x)\dd G(x)=0$. Thus \eqref{def_vague} does not hold. However, it is interesting to note that  there is a convergence of Fourier--Stieltjes transforms. Indeed, for any $t\in\R$
\begin{eqnarray*}
	\int_{\R} e^{itx} \dd G_n(x)=\bigl(1-e^{it/n^2}\bigr)\cdot n= -\dfrac{it}{n}\,(1+o(1))\to 0=\int_{\R} e^{itx} \dd G(x),\quad n\to\infty.\quad \Box\\
\end{eqnarray*}

The next example shows that the use of the sebsequences in the definition of basic convergence is essential.

\textbf{Example 3.} For any $n\in\N$ we set $k_n\in\N\cup\{0\}$ satisfying $2^{k_n}\leqslant n<2^{k_n+1}$.  We define
 \begin{eqnarray*}
 	G_n(x)\defeq\id_{a_n}(x)-\id_{b_n}(x),\quad x\in\R, \quad \text{where}\quad a_n\defeq\dfrac{n-2^{k_n}}{2^{k_n}}\quad\text{and}\quad b_n\defeq\dfrac{n+1-2^{k_n}}{2^{k_n}},\quad n\in\N.
\end{eqnarray*}
It is seen that the interval $[a_n, b_n]$ is vanishing ($b_n-a_n=2^{-k_n}\to 0$) and shifting over $[0,1]$ as $n\to\infty$. 
Let $h:\R\to\R$ be a bounded continuous function. Due to the uniform continuity of $h$ on $[0,1]$, we have that
\begin{eqnarray*}
	\int_{\R} h(x)\dd G_n(x)=h(a_n)-h(b_n)\to 0,\quad n\to\infty. 
\end{eqnarray*}
So $(G_n)_{n\in\N}$  weakly converges to $G(x)\defeq0$ for all $x\in\R$. Then, by the comments above, $(G_n)_{n\in\N}$ basically converges to $G$ that can be also checked directly by definition. However, for any $x_0,x_1, x_2\in[0,1)$ there is no limit either for $G_n(x_0)$ or $G_n(x_2)-G_n(x_1)$ as $n\to\infty$, because $G_n(x)$ takes an infinite number of times each of the values  $1$ or $0$, when $x\in[a_n,b_n)$ or not correspondingly. Note that, due to the weak convergence of $(G_n)_{n\in\N}$, there is a convergence of Fourier--Stieltjes transforms:  $\int_{\R} e^{itx} \dd G_n(x)\to \int_{\R} e^{itx} \dd G(x)=0$ as $n\to\infty$ for every $t\in\R$. \quad $\Box$\\

The following example shows that  it is important to use the differences of values of the functions at points $x_1$ and $x_2$ in the definition of the basic convergence in order to stay within $\BV$.
  
\textbf{Example 4.} For every $n\in\N$ we define $G_n(x)\defeq 1+\tfrac{x}{n}$ for $x\in[-n,n]$, $G_n(x)=0$ for $x<-n$, and  $G_n(x)=1$ for $x > n$.  So $G_n$ are  non-decreasing continuous functions, and $G_n(+\infty)=\|G_n\|=2$, $n\in\N$. We see that  $G_n(x)\to 1$ as  $n\to\infty$ for any $x\in\R$. However, an identical 1  doesn't belong to $\BV$ (it must be $0$ at $-\infty$).  At the same time for any real $x_1$ and $x_2$ we have $G_n(x_2)-G_n(x_1)\to 0$, $n\to\infty$, and we conclude that $(G_n)_{n\in\N}$ basically converges to the function $G(x)\defeq 0$ for all $x\in\R$, which is from $\BV$. Of course, $(G_n)_{n\in\N}$ doesn't weakly converge to $G$ here, because $\int_{\R} \dd G_n(x)=2\nrightarrow \int_{\R} \dd G(x)=0$, $n\to\infty$.

Note that for $t\ne0$
\begin{eqnarray*}
	\int_{\R} e^{itx} \dd G_n(x)= \dfrac{1}{n}\int_{-n}^n e^{itx} \dd x= \dfrac{e^{itn}-e^{-itn}}{it n}\to 0,\quad n\to \infty,
\end{eqnarray*}
and
\begin{eqnarray*}
	\int_{\R} e^{itx} \dd G_n(x)\biggr|_{t=0}= \int_{\R}  \dd G_n(x) =G_n(+\infty)=2,\quad n\in\N.
\end{eqnarray*}
Thus the Fourier--Stieltjes transforms of $G_n$, $n\in\N$, pointwisely converge to the Fourier--Stieltjes transform of $G$ (i.e. to identical~0) for almost all $t\in\R$.\quad $\Box$\\

We now consider a general question about the relationship between the basic convergence of functions from $\BV$ and the convergence of their Fourier--Stieltjes transforms. We are not pretended to full studing of this question here, and we present only those assertions that will be used below in the main results of the article.

Let $(G_n)_{n\in\N}$ be a  sequence of functions from $\BV$. Let us define the corresponding sequence of Fourier--Stieltjes integrals:
\begin{eqnarray*}
	g_n(t)=\int_{\R} e^{itx} \dd G_n(x),\quad t\in \R,\quad  n\in\N.
\end{eqnarray*}
The results below in fact show that, under the rather weak and natural  assumptions,  the pointwise convergence of $g_n$ implies the basic convergence of $G_n$ as $n\to\infty$.

We will use the following assumption:
\begin{eqnarray}\label{assum_UBV}
	\varlimsup_{n\to\infty} \|G_n\|=B<\infty.
\end{eqnarray}

\begin{Theorem}\label{th_FStoFBV}
	Let $(G_n)_{n\in\N}$ satisfy \eqref{assum_UBV}. 	Suppose that $g_n(t)\to g(t)$, $n\to\infty$, for almost all $t\in\R$ with some function $g:\R\to\Complex$. Then there exists a function $G\in\BV$  such that $\|G\|\leqslant B$ and the equality
	\begin{eqnarray}\label{th_FStoFBV_conc}
		g(t)=\int_{\R} e^{itx} \dd G(x)
	\end{eqnarray}
	holds for almost all $t\in\R$ including  all continuity points of the function $g$. The function $G$ is uniquely determined in the class $\BV$, and $G_n\cconv G$, $n\to\infty$. If also $g_n(0)\to g(0)$, $n\to\infty$, and $g$ is continuous at $t=0$, then additionally $G_n(+\infty)\to G(+\infty)$, $n\to\infty$.
\end{Theorem}

We  are not aware of any results with such assertion. There are some close remarks in \cite{Bochner} and \cite{Phillips}. It is seen that this theorem complements and partially generalizes the well known Levy's continuity theorem, which was stated for sequences of probability distribution functions. 

Suppose that the sequence $(G_n)_{n\in\N}$  weakly converges to a function $G$ from $\BV$ with Fourier--Stieltjes transform $g$. Then \eqref{assum_UBV} is satisfied (see Proposition 1.4.4. in \cite{Bogachev}, p. 22) and $g_n(t)\to g(t)$, $n\to\infty$, for every $t\in\R$. According to the theorem, we have the basic  convergence $G_n\cconv G$ and also $G_n(+\infty)\to G(+\infty)$, $n\to\infty$. Thus we showed that the weak convergence implies the basic convergence.

We next formulate the analog of Theorem \ref{th_FStoFBV} with using of the decompositions
\begin{eqnarray}\label{eq_fbvn_decom}
	G_n(x)=G_n^+(x)-G_n^-(x),\quad x\in\R,\quad n\in\N,
\end{eqnarray}
where $G_n^+$ and $G_n^-$ are non-decreasing functions from $\BV$. Here we will assume:
\begin{eqnarray}\label{assum_UBVNeg}
	\varlimsup_{n\to\infty} G_n^-(+\infty)=M<\infty,
\end{eqnarray}
that  can sometimes be  more convenient for checking than \eqref{assum_UBV}.

\begin{Proposition}\label{pr_FS_0}
	Let $(G_n)_{n\in\N}$ satisfy \eqref{assum_UBVNeg} for some decompositions \eqref{eq_fbvn_decom}. Suppose that $g_n(0)\to c\in\R$, $n\to\infty$. Then  \eqref{assum_UBV} holds with some $B\leqslant c+2M$.
\end{Proposition}
Thus if there is a convergence of $g_n$, $n\in\N$, at $t=0$, then for some \eqref{eq_fbvn_decom} assumptions \eqref{assum_UBV}  and \eqref{assum_UBVNeg}  are equivalent. So we come to the following assertion. 

\begin{Theorem}\label{th_FStoFBV_0}
	Let $(G_n)_{n\in\N}$ satisfy \eqref{assum_UBVNeg} for some decompositions \eqref{eq_fbvn_decom}. Suppose that $g_n(t)\to g(t)$, $n\to\infty$, for almost all $t\in\R$ including $t=0$ with some function $g:\R\to\Complex$. Then $g(0)\in\R$, the condition  \eqref{assum_UBV} is satisfied with $B\leqslant g(0)+2M$, and the assertions of Theorem \ref{th_FStoFBV} hold.  If also $g$ is continuous at $t=0$, then additionally $G_n(+\infty)\to G(+\infty)$, $n\to\infty$.
\end{Theorem}

All these results will be proved in Section 4.

\section{Main results}

Let $(F_n)_{n\in\N}$ be a sequence of quasi-infinitely divisible distribution function with corresponding sequence of characteristic function $(f_n)_{n\in\N}$. Let every $f_n$ admit the representation
\begin{eqnarray}\label{def_fn}
	f_n(t)=\exp\Biggl\{it \gamma_n+\int_{\R} \Bigl(e^{itx} -1 -\tfrac{it}{\tau} \sin(\tau x)\Bigr)\tfrac{1+x^2 }{x^2}  \dd G_n(x)\Biggr\},\quad t\in\R,\quad n\in\N,
\end{eqnarray}
where  $\gamma_n\in\R$, $G_n\in \BV$, $n\in\N$, and $\tau>0$ is a fixed number. We are interested in criteria of the weak convergence of $(F_n)_{n\in\N}$ in terms of the spectral pairs $(\gamma_n, G_n)$, $n\in\N$.

Assertions of the following Theorems \ref{th_LPS1} and \ref{th_LPS2} were obtained by Lindner, Pan, and Sato in \cite{LindPanSato} (where the results were presented in another form). 

\begin{Theorem}\label{th_LPS1}
	If $\gamma_n\to \gamma$ and  $G_n\xrightarrow{w} G$, $n\to\infty$, with some $\gamma\in\R$ and $G\in\BV$, then $(\gamma, G)$ is the spectral pair for some quasi-infinitely divisible distribution function $F$, and $(F_n)_{n\in\N}$ weakly converges to $F$.
\end{Theorem}

We next use the decompositions
\begin{eqnarray}\label{eq_Gn_decom}
	G_n(x)=G_n^+(x)-G_n^-(x),\quad x\in\R,\quad n\in\N.
\end{eqnarray}
where $G_n^+$ and $G_n^-$ are non-decreasing functions from $\BV$. There exists an important  way of choosing $G_n^+$ and $G_n^-$. Let $\mu_{G_n}$ be the signed measure that is generated by $G_n$ for every $n\in\N$, i.e. such that $\mu_{G_n}((a,b])=G_n(b)-G_n(a)$ for all $a,b\in\R$, $a\leqslant b$, $n\in\N$. Every measure $\mu_{G_n}$ is uniquely represented by \textit{the Hahn--Jordan
decomposition} $\mu_{G_n}=\mu_{G_n}^+ -\mu_{G_n}^-$, where $\mu_{G_n}^+ $ and $\mu_{G_n}^-$ are non-negative finite measures concentrated on some disjoint sets (see \cite{Bogachev} p. 3). So we can choose
\begin{eqnarray}\label{def_HahnJordan}
	G_n^+(x)= \mu_{G_n}^+((-\infty,x]),\quad\text{and} \quad G_n^-(x)= \mu_{G_n}^-((-\infty,x]),\quad x\in\R,\quad n\in\N.
\end{eqnarray}
In this case we will  have \eqref{eq_Gn_decom} and additionally that $|G_n|(x)=G_n^+(x)+G_n^-(x)$, $x\in\R$, $n\in\N$.

\begin{Theorem}\label{th_LPS2}
	Let $F$ be a distribution function and $(F_n)_{n\in\N}$ weakly converges to $F$. Suppose that $G_n^+$ and $G_n^-$ from \eqref{eq_Gn_decom} are defined according to the Hahn--Jordan
	decomposition by \eqref{def_HahnJordan} for every $n\in\N$. Suppose that the sequence $(G^-_n)_{n\in\N}$ satisfies the assumptions 
	\begin{eqnarray*}
		\sup_{n\in\N} \|G^-_n\|<\infty\quad\text{and}\quad \lim_{r\to\infty}\sup_{n\in\N} \bigl(1-|G_n^-|(r)+|G_n^-|(-r)\bigr)=0
	\end{eqnarray*}
	$($uniform boundedness in variation and tightness, correspondingly\,$)$. Then $F$ is quasi-infinitely divisible  with some spectral pair $(\gamma, G)$. Moreover, $\gamma_n\to \gamma$ and $G_n\xrightarrow{w} G$, $n\to\infty$.
\end{Theorem}

Theorems \ref{th_LPS1} and \ref{th_LPS2} connect the weak convergence of quasi-infinitely divisible distribution functions with the weak convergence of their spectral functions.  We are interested in analogs of these theorems but with the basic convergence of the spectral functions.

We will use the following assumption:
\begin{eqnarray}\label{assum_UNB_Gn}
	\varlimsup_{n\to\infty} \|G_n\|=B<\infty.
\end{eqnarray}

\begin{Theorem}\label{th_FnconvtoF}
	Suppose that $(F_n)_{n\in\N}$  satisfies \eqref{assum_UNB_Gn} with some $B\geqslant 0$. Let  $(F_n)_{n\in\N}$ weakly converge to a distribution function $F$. Then $F$ is quasi-infinitely divisible with some spectral pair $(\gamma, G)$, where $\gamma\in\R$ and $G\in\BV$ with $\|G\|\leqslant B$.   Moreover,  $\gamma_n\to \gamma$, $G_n\cconv G$, and $G_n(+\infty)\to G(+\infty)$, $n\to\infty$.
\end{Theorem}

The next theorem is an analog of this one, but with the assumption
\begin{eqnarray}\label{assum_GnNeg}
	\varlimsup_{n\to\infty} G_n^-(+\infty)=M<\infty,
\end{eqnarray}
on decompositions \eqref{eq_Gn_decom} for $G_n$, $n\in\N$.  If we choose  $G_n^+$ and $G_n^-$ according to the Hahn--Jordan
decomposition by \eqref{def_HahnJordan} for every $n\in\N$, then \eqref{assum_GnNeg} is weaker than \eqref{assum_UNB_Gn}. Also observe that \eqref{assum_GnNeg} is satisfied, when we deal with non-decreasing functions $G_n$, $n\in\N$. It should be noted, however, that it is not required in the theorems and corollaries below that  $G_n^+$ and $G_n^-$  in  \eqref{eq_Gn_decom} must be choosen according to the Hahn--Jordan decomposition.

\begin{Theorem}\label{th_FnconvtoFminus}
	 Suppose that $(F_n)_{n\in\N}$  satisfies \eqref{assum_GnNeg} with some $M\geqslant 0$ and for some decompositions \eqref{eq_Gn_decom}.  Let  $(F_n)_{n\in\N}$ weakly converge to a distribution function $F$. Then \eqref{assum_UNB_Gn} holds for some $B\geqslant 0$ and all assertions of Theorem \ref{th_FnconvtoF} are true. Also we have that  $B\leqslant G(+\infty)+2M$. 
\end{Theorem}

Theorems \ref{th_FnconvtoF} and \ref{th_FnconvtoFminus} yield necessary conditions for the weak convergence within the class of  quasi-infinitely divisible distribution functions under the assumption \eqref{assum_UNB_Gn} or \eqref{assum_GnNeg}. 

\begin{Corollary}\label{co_FnconvF_ness}
	Suppose that $(F_n)_{n\in\N}$  satisfies \eqref{assum_UNB_Gn} or \eqref{assum_GnNeg} for some decompositions \eqref{eq_Gn_decom}. Let $F$ be a quasi-infinitely divisible distribution function $F$ with spectral pair $(\gamma, G)$, where $\gamma\in\R$ and $G\in\BV$. If the sequence $(F_n)_{n\in\N}$ weakly converges to $F$, then $\gamma_n\to \gamma$, $G_n\cconv G$, and $G_n(+\infty)\to G(+\infty)$, $n\to\infty$.
\end{Corollary}

Also Theorems \ref{th_FnconvtoF} and \ref{th_FnconvtoFminus} state  sufficient conditions for membership of the class of quasi-infinitely divisible distribution functions. 

\begin{Corollary}\label{co_suff}
	A distribution function $F$ is quasi-infinitely divisible if it is a weak limit of a sequence $(F_n)_{n\in\N}$ of quasi-infinitely divisible distribution functions $($with characteristic functions \eqref{def_fn}$)$, which satisfies \eqref{assum_UNB_Gn} or \eqref{assum_GnNeg} for some decompositions \eqref{eq_Gn_decom}.
\end{Corollary}

Note that this corollary is a stronger version of the same assertion in Theorem \ref{th_LPS2}, because we don't assume the tightness for $(G_n^-)_{n\in\N}$ and we don't require the use of  the Hahn--Jordan decomposition.

It is known (see \cite{LindPanSato} p.  17) that a weak limit of quasi-infinitely divisible distribution function is not necessarily quasi-infinitely divisible. Hence assumptions \eqref{assum_UNB_Gn} or \eqref{assum_GnNeg} can not be simply omitted in Corollary~\ref{co_suff}.
However, it seems that they can be done weaker (see \cite{LindPanSato}  Example 4.4).

We will use a notion of relative compactness for $(F_n)_{n\in\N}$ in the next theorem. Recall that  $(F_n)_{n\in\N}$ is said to be \textit{relatively compact} if every its subsequence contains a further subsequence that weakly converges to a distribution function. It is clear that a weakly convergent sequence of distribution functions is relatively compact.
In general, the property of relative compactness is not difficult for checking due to the Prokhorov's theorem and various probability inequalities. Also some criteria of relative compactness are known for particular important sequences of distribution functions (for example, see \cite{KhartovComp1}, \cite{KhartovComp2}, \cite{Khartov}, and the references given there).

\begin{Theorem}\label{th_FnconvtoF_suff}
	Suppose that $(F_n)_{n\in\N}$  satisfies \eqref{assum_UNB_Gn}. If $(F_n)_{n\in\N}$ is relatively compact and $\gamma_n\to \gamma$,  $G_n\cconv G$, $n\to\infty$, with some $\gamma\in\R$ and $G\in\BV$, then $(\gamma, G)$ is the spectral pair for a quasi-infinitely divisible distribution function $F$ and the sequence $(F_n)_{n\in\N}$ weakly converges to $F$.
\end{Theorem}

This theorem yields sufficient conditions for the weak convergence within the class of  quasi-infinitely divisible distribution functions under the assumption \eqref{assum_UNB_Gn}.

\begin{Corollary}\label{co_FnconvF_suff}
	Suppose that $(F_n)_{n\in\N}$  satisfies \eqref{assum_UNB_Gn}. Let $F$ be a quasi-infinitely divisible distribution function $F$ with spectral pair $(\gamma, G)$, where $\gamma\in\R$ and $G\in\BV$.  If $(F_n)_{n\in\N}$ is relatively compact and  $\gamma_n\to \gamma$,  $G_n\cconv G$, $n\to\infty$, then $(F_n)_{n\in\N}$ weakly converges to $F$.
\end{Corollary}

Corollaries \ref{co_FnconvF_ness} and \ref{co_FnconvF_suff} directly yield the following criterion.

\begin{Theorem}\label{th_FnconvtoF_crit}
	Suppose that $(F_n)_{n\in\N}$  satisfies \eqref{assum_UNB_Gn}. Let $F$ be a quasi-infinitely divisible distribution function $F$ with spectral pair $(\gamma, G)$, where $\gamma\in\R$ and $G\in\BV$. The sequence $(F_n)_{n\in\N}$ weakly converges to $F$ if and only if  $(F_n)_{n\in\N}$ is relatively compact and  $\gamma_n\to \gamma$,  $G_n\cconv G$, $n\to\infty$. Moreover, the convergence  $G_n(+\infty)\to G(+\infty)$, $n\to\infty$, can be added to the necessary conditions.
\end{Theorem}

We now formulate the analogs of Theorems \ref{th_FnconvtoF_suff} and \ref{th_FnconvtoF_crit}, and of  Corollary \ref{co_FnconvF_suff} under the assumption \eqref{assum_GnNeg}. They are directly stated due to the following simple note.  

Suppose that a sequence $(G_n)_{n\in\N}$ from $\BV$ satisfies \eqref{assum_GnNeg} for some decompositions~\eqref{eq_Gn_decom}. If $\varlimsup_{n\to\infty} G_n(+\infty)$  is finite, then \eqref{assum_UNB_Gn} holds. Indeed, according to \eqref{eq_Gn_decom}, it follows from the inequalities
\begin{eqnarray*}\label{ineq_varGn}
	\|G_n\|\leqslant \|G_n^+\|+\|G_n^-\|= G_n^+(+\infty)+G_n^-(+\infty)=G_n(+\infty)+2 G_n^-(+\infty),\quad n\in\N.
\end{eqnarray*}
 
So we obtain the following results.

\begin{Theorem}\label{th_FnconvtoFminus_suff}
	Suppose that $(F_n)_{n\in\N}$  satisfies \eqref{assum_GnNeg} for some decompositions \eqref{eq_Gn_decom}. If $(F_n)_{n\in\N}$ is relatively compact and $\gamma_n\to \gamma$,  $G_n\cconv G$, and $G_n(+\infty)\to G(+\infty)$, $n\to\infty$, with some $\gamma\in\R$ and $G\in\BV$, then all assertions of Theorem \ref{th_FnconvtoF_suff} hold.
\end{Theorem}

So  Theorems \ref{th_FnconvtoF_suff} and \ref{th_FnconvtoFminus_suff} complement Theorem \ref{th_LPS1}: we use weaker convergence for the spectral functions $(G_n)_{n\in\N}$, but we additionally assume the relative compactness of $(F_n)_{n\in\N}$.

\begin{Corollary}\label{co_FnconvFminus_suff}
	Suppose that $(F_n)_{n\in\N}$  satisfies \eqref{assum_GnNeg} for some decompositions \eqref{eq_Gn_decom}. Let $F$ be a quasi-infinitely divisible distribution function $F$ with the spectral pair $(\gamma, G)$, where $\gamma\in\R$ and $G\in\BV$.  If $(F_n)_{n\in\N}$ is relatively compact and  $\gamma_n\to \gamma$,  $G_n\cconv G$, and $G_n(+\infty)\to G(+\infty)$, $n\to\infty$, then $(F_n)_{n\in\N}$ weakly converges to $F$.
\end{Corollary}

Corollaries \ref{co_FnconvF_ness} and \ref{co_FnconvFminus_suff} directly yield the following criterion.
\begin{Theorem}\label{th_FnconvtoFminus_crit}
	Suppose that $(F_n)_{n\in\N}$  satisfies \eqref{assum_GnNeg} for some decompositions \eqref{eq_Gn_decom}. Let $F$ be a quasi-infinitely divisible distribution function $F$ with spectral pair $(\gamma, G)$, where $\gamma\in\R$ and $G\in\BV$. The sequence $(F_n)_{n\in\N}$ weakly converges to $F$ if and only if  $(F_n)_{n\in\N}$ is relatively compact and  $\gamma_n\to \gamma$,  $G_n\cconv G$, and $G_n(+\infty)\to G(+\infty)$, $n\to\infty$.
\end{Theorem}

On account of comments before Example 1 in Section 2, Theorems \ref{th_FnconvtoF_crit} and \ref{th_FnconvtoFminus_crit} complement  well known similar results for the weak convergence of infinitely divisible distribution functions (see \cite{GnedKolm} p. 87). 

\section{Proofs}

\textbf{Proof of Theorem \ref{th_FStoFBV}.}\quad First, observe that the function $g$ is measurable, because it is a almost everywhere limit of continuous (hence measurable) functions $g_n$, $n\in\N$. So we have
\begin{eqnarray}\label{conc_convrho}
	\int_{\R} g_n(t)\rho(t) \dd t \to \int_{\R} g(t)\rho(t) \dd t,\quad n\to\infty,
\end{eqnarray}
for any function $\rho\in L_1(\R)$. Indeed, due to \eqref{assum_UBV}, there exists a constant $B_0>0$ such that $|g_n(t)|\leqslant \|G_n\|\leqslant B_0$ for all $n\in\N$, and convergence \eqref{conc_convrho} holds by the Lebesgue dominated convergence theorem.

Let us define the function
\begin{eqnarray*}
	\varphi(x)\defeq \int_{\R} e^{itx} \rho(t) \dd t,\quad x\in\R.
\end{eqnarray*}
Observe that for every $n\in\N$ we have
\begin{eqnarray}\label{eq_int_fsn_rho}
	\int_{\R} g_n(t)\rho(t) \dd t = 	\int_{\R}\biggl(\int_{\R} e^{itx} \dd G_n(x)\biggr) \rho(t) \dd t=\int_{\R} \biggl(\,\, \int_{\R} e^{itx}\rho(t) \dd t\biggr)  \dd G_n(x)=\int_{\R} \varphi(x)  \dd G_n(x).
\end{eqnarray}
Let us consider the last integral.  Due to \eqref{ineq_fbv_var} and \eqref{assum_UBV}, by  Helly's first theorem (see \cite{Natanson}, pp. 222 and 240), there exists a subsequence $(G_{n_k})_{k\in\N}$ in $(G_n)_{n\in\N}$ and a function of bounded variation $G_*:\R\to\R$ such that $G_{n_k}(x)\to G_*(x)$ as $k\to\infty$ for all $x\in\R$. Note that, in general,  $G_*$ may not be right-continuous (see Example 3). But  $\varphi$ is bounded and continuous on $\R$ and hence there exists the Riemann--Stieltjes  integral $\int_{\R} \varphi(x)\dd G_*(x)$. Also the (Lebesgue--Stieltjes) integrals $\int_{\R} \varphi(x) \dd G_n(x)$ coincide with the corresponding Riemann--Stieltjes integrals. Next, due to the well known fact that $\varphi(x)\to 0$ as $x\to\pm \infty$, by Helly's second theorem (see \cite{Natanson}, p. 240), we have the following convergence for the Riemann--Stieltjes integrals:
\begin{eqnarray*}
	\int_{\R} \varphi(x)  \dd G_{n_k}(x)\to \int_{\R} \varphi(x)  \dd G_*(x),\quad k\to\infty.
\end{eqnarray*}
Let us define $G(x)\defeq G_*(x+0)-G_*(-\infty)$, $x\in\R$ (note that $G_*(-\infty)\ne 0$ in general, see Example 4). So $G$ is right-continuous on $\R$ and $G(-\infty)=0$, i.e. $G\in\BV$. Since  $G(x)$ equals $G_*(x)-G_*(-\infty)$ for all  $x\in\R$ except at most countable set, due to the continuity of $\varphi$, we have 
\begin{eqnarray*}
	\int_{\R}\varphi(x)  \dd G_*(x)= \int_{\R} \varphi(x)  \dd G(x),
\end{eqnarray*}
where the integral in the right-hand side can be considered as Lebesgue--Stieltjes integral. Thus we have the following convergence with the Lebesgue--Stieltjes  integrals:
\begin{eqnarray*}
	\int_{\R} \varphi(x)  \dd G_{n_k}(x)\to \int_{\R} \varphi(x)  \dd G(x),\quad k\to\infty.
\end{eqnarray*}
The integral in the right-hand side admits the following representation analogously to \eqref{eq_int_fsn_rho}:
\begin{eqnarray*}
	\int_{\R}\varphi(x)  \dd G(x)=\int_{\R} \biggl(\int_{\R} e^{itx} \dd G(x)  \biggr) \rho(t) \dd t.
\end{eqnarray*}
Due to \eqref{conc_convrho} and \eqref{eq_int_fsn_rho}, we also have
\begin{eqnarray*}
	\int_{\R} \varphi(x)  \dd G_{n_k}(x)\to \int_{\R} g(t)\rho(t) \dd t,\quad k\to\infty.
\end{eqnarray*}
Thus we obtain
\begin{eqnarray}\label{eq_fs_fbv_rho}
	\int_{\R} g(t)\rho(t) \dd t=\int_{\R} \biggl(\int_{\R} e^{itx} \dd G(x)  \biggr) \rho(t) \dd t
\end{eqnarray}
for any function  $\rho\in L_1(\R)$. This implies that
\begin{eqnarray}\label{eq_fs_int}
	g(t)=\int_{\R} e^{itx} \dd G(x)\quad\text{for almost every}\,\,\, t\in\R.
\end{eqnarray}
Indeed,  conversely, suppose that there exists a bounded set $E$ of non-zero Lebesgue measure such that $\Delta(t)\defeq g(t)-\int_{\R} e^{itx} \dd G(x)\ne 0$, $t\in E$. Let us introduce the sets 
\begin{eqnarray*}
	&&E_1\defeq\{t\in E: \Realpart \Delta(t)>0\},\quad E_2\defeq \{t\in E: \Realpart \Delta(t)<0\},\\
	&&E_3\defeq\{t\in E: \Imagpart \Delta(t)>0\},\quad E_4\defeq\{t\in E: \Imagpart \Delta(t)<0\}.
\end{eqnarray*}
It easily seen that $E=E_1\cup E_2\cup E_3\cup E_4$. Hence at least one $E_j$ has non-zero Lebesgue measure. We denote any such set by $E_*$. Next, according to the property of strict positivity of integral, we obtain
\begin{eqnarray*}
	\biggl|\int_{E_*} \Delta (t) \dd t \biggr|\geqslant \biggl|\int_{E_*} \Realpart\Delta (t) \dd t \biggr|=\int_{E_*} \bigl|\Realpart\Delta (t)\bigr| \dd t>0, \quad\text{for}\quad E_*=E_1\,\text{ or }\, E_*=E_2,
\end{eqnarray*}
and
\begin{eqnarray*}
	\biggl|\int_{E_*} \Delta (t) \dd t \biggr|\geqslant \biggl|\int_{E_*} \Imagpart\Delta (t) \dd t \biggr|=\int_{E_*} \bigl|\Imagpart\Delta (t)\bigr| \dd t>0,\quad\text{for}\quad E_*=E_3\,\text{ or }\, E_*=E_4.
\end{eqnarray*}
Thus we have
\begin{eqnarray*}
	\Biggl| \int_{E_*} g(t) \dd t-\int_{E_*} \biggl(\int_{\R} e^{itx} \dd G(x)  \biggr)  \dd t \Biggr|=\biggl|\int_{E_*} \Delta (t) \dd t \biggr|>0.
\end{eqnarray*}
This contradicts  \eqref{eq_fs_fbv_rho} when we choose $\rho$ as follows:  $\rho(t)=1$, $t\in E_*$, and $\rho(t)=0$, $t\notin E_*$. It is valid since $\rho\in L_1(\R)$ due to the boundedness $E_*\subset E$. Thus \eqref{eq_fs_int} is true.

Let us show that \eqref{eq_fs_int} holds for every continuity point of the function $g$. Let $\aeset$ be the set of all $t\in\R$ for which  \eqref{eq_fs_int} holds. Hence the Lebesgue measure of $\R\setminus\! \aeset$ equals zero.  Let $g$ be continuous at the fixed point $t_0$.  So we can choose  $t_m\in\aeset$, $m\in\N$ such that $t_m\to t_0$, $m\to\infty$. Then 
$g(t_m)\to g(t_0)$, $m\to\infty$, and  at the same time
\begin{eqnarray*}
	g(t_m)=\int_{\R} e^{it_m x} \dd G(x)\to \int_{\R} e^{it_0 x} \dd G(x),\quad m\to\infty.
\end{eqnarray*}
due to continuity of the function $t\mapsto \int_{\R} e^{it x} \dd G(x)$ on $\R$. Thus we have $g(t_0)=\int_{\R} e^{it_0 x} \dd G(x)$.

According to equality \eqref{eq_fs_int}, the function $g$  almost everywhere coincides with the continuous function $t\mapsto\int_{\R} e^{it x} \dd G(x)$, $t\in\R$. So the latter  function is uniquely determined by $g$ within the class of all continuous complex-valued functions on $\R$. Next, it is well known that $t\mapsto\int_{\R} e^{it x} \dd G(x)$, $t\in\R$, uniquely determines $G$ within the class $\BV$. Therefore $g$ uniquely determines $G$ in the class $\BV$.

Let us return to the sequence $(G_{n_k})_{k\in\N}$. From the above we know that $G_{n_k}(x)\to G_*(x)$ for all $x\in\R$, and $G(x)=G_*(x)-G_*(-\infty)$ for all  $x\in\R$ except at most a countable set $D$ where $G_*$ is not right-continuous. Then for all $x_1,x_2\in\R\setminus\! D$ we have 
\begin{eqnarray}\label{conc_conv_fbvsubseq}
	G_{n_k}(x_2)-G_{n_k}(x_1)\to \bigl(G(x_2)+G_*(-\infty)\bigr) - \bigl(G(x_1)+G_*(-\infty)\bigr)=G(x_2)-G(x_1),\quad k\to\infty.
\end{eqnarray}

Let  $(G_{m_l})_{l\in\N}$ be an arbitrary subsequence of $(G_n)_{n\in\N}$. Analogously to the above,  there exists a further subsequence $(G_{m'_k})_{k\in\N}$ in $(G_{m_k})_{k\in\N}$, which pointwise converges to some function of bounded variation $H_*:\R\to\R$, i.e. $G_{m'_k}(x)\to H_*(x)$, $k\to\infty$, for all $x\in\R$. Defining $H(x)\defeq H_*(x+0)-H_*(-\infty)$, $x\in\R$, we as before will obtain  $g(t)=\int_{\R} e^{itx}\dd H(x)$ for almost all $t\in\R$, with $H\in\BV$. Since $G$ is a unique function within $\BV$, which represents  $g$ by \eqref{th_FStoFBV_conc}, we have $H(x)=G(x)$, $x\in\R$. We also have
\begin{eqnarray*}
	G_{m'_k}(x_2)-G_{m'_k}(x_1)\to G(x_2)-G(x_1),\quad k\to\infty,
\end{eqnarray*}
for all $x_1,x_2\in\R$ except at most countable set $D'$ where $H_*$ is not right-continuous (in general $D'\ne D$).  So we proved that $G_n\cconv G$, $n\to\infty$.

Let us consider the numbers  $g_n(0)=\int_{\R}\dd G_n(x)=G_n(+\infty)$, $n\in\N$. If we suppose that $g$ is continuous at $t=0$, then, by the above remarks, we will have $g(0)=\int_{\R}\dd G(x)=G(+\infty)$. Therefore, assuming to hold $g_n(0)\to g(0)$, $n\to\infty$, we will obtain $G_n(+\infty)\to G(+\infty)$, $n\to\infty$.

It remains to prove that $\|G\|\leqslant B$. On the contrary, suppose that this is false. Then we can find $y_0, y_1,\ldots, y_N\in\R$ such that
\begin{eqnarray}\label{assum_sumvariation}
	B<\sum\limits_{j=1}^{N}\bigl|G(y_j)-G(y_{j-1})\bigr|\leqslant\|G\|.
\end{eqnarray}
Let us take our sequence $(G_{n_k})_{k\in\N}$ and the set $D$, which is at most countable. Since $G$ is right-continuous and the set  $\R\setminus\! D$ is dense, we can assume that $y_0$, $y_1$, \ldots, $y_N$ are choosen from $\R\setminus\!D$. Next, due to the convergence \eqref{conc_conv_fbvsubseq} and assumption \eqref{assum_UBV}, we have 
\begin{eqnarray*}
	\sum\limits_{j=1}^{N}|G(y_j)-G(y_{j-1})|=\lim\limits_{k\to\infty}\sum\limits_{j=1}^{N}|G_{n_k}(y_j)-G_{n_k}(y_{j-1})|\leqslant \varlimsup_{n\to\infty} \|G_{n}\|\leqslant B,
\end{eqnarray*}
which contradicts \eqref{assum_sumvariation}.\quad $\Box$\\

\textbf{Proof of Proposition \ref{pr_FS_0}.}\quad  By the assumption $g_n(0)\to c\in\R$, $n\to\infty$. Since $g_n(0)=\int_{\R}\dd G_n(x)=G_n(+\infty)$, $n\in\N$, we have the convergence $G_n(+\infty)\to c$, $n\to\infty$. Let us consider decompositions \eqref{eq_fbvn_decom}. We have $G_n(+\infty)=G^+_n(+\infty)-G^-_n(+\infty)$, $n\in\N$. Also observe that
\begin{eqnarray*}
	\|G_n\|\leqslant \|G^+_n\| +\|G^-_n\|=G^+_n(+\infty)+G^-_n(+\infty)=G_n(+\infty)+2G^-_n(+\infty),\quad n\in\N. 
\end{eqnarray*}
Therefore
\begin{eqnarray*}
	B=\varlimsup_{n\to\infty} \|G_n\|\leqslant \lim_{n\to\infty}G_n(+\infty)+2\varlimsup_{n\to\infty}G^-_n(+\infty)
	=c+2M.
\end{eqnarray*}
Thus we have \eqref{assum_UBV} with $B\leqslant g(0)+2M$. \quad $\Box$\\

\textbf{Proof of Theorem \ref{th_FStoFBV_0}.}\quad  By the assumption $g_n(0)\to g(0)$, $n\to\infty$. Since $g_n(0)=\int_{\R}\dd G_n(x)=G_n(+\infty)$, $n\in\N$, we have the convergence $G_n(+\infty)\to g(0)$, $n\to\infty$. So the sequence $G_n(+\infty)\in\R$, $n\in\N$,  has a finite limit $g(0)$ that must be real. According to Proposition \ref{pr_FS_0}, condition \eqref{assum_UBV} holds with some $B\leqslant g(0)+2M$. Using Theorem \ref{th_FStoFBV}, we get all its assertions. So $g(t)=\int_{\R} e^{itx} \dd G(x)$ holds for some $G\in\BV$ and for all $t\in\R$ that are continuity points of the function $g$. Under the assumption, $g$ is continuous at $t=0$, and we have $g(0)=\int_{\R}\dd G(x)=G(+\infty)$. Since $G_n(+\infty)\to g(0)$, $n\to\infty$, we obtain that $G_n(+\infty)\to G(+\infty)$, $n\to\infty$.\quad $\Box$\\

We need the following lemma for proving Theorem \ref{th_FnconvtoF}. 
\begin{Lemma}\label{lm_UVW} For any $t\in\R$ and $\tau>0$ the following representations hold
	\begin{eqnarray}
		&&e^{itx} -1 -\tfrac{it}{\tau} \sin(\tau x)=\int_{A_{t,\tau}} e^{isx} \dd U_{t,\tau}(s),\qquad \Bigl(e^{itx} -1 -\tfrac{it}{\tau} \sin(\tau x)\Bigr)\dfrac{1}{x^2}=\int_{A_{t,\tau}} e^{isx} \dd V_{t,\tau}(s),\label{lm_UVW_intUV}\\
		&&\qquad\qquad\qquad\Bigl(e^{itx} -1 -\tfrac{it}{\tau} \sin(\tau x)\Bigr)\tfrac{1+x^2}{x^2}=\int_{A_{t,\tau}} e^{isx} \dd W_{t,\tau}(s)\label{lm_UVW_intW},\quad x\in\R,
	\end{eqnarray}
	where $A_{t,\tau}\defeq \bigl\{s\in\R: |s|\leqslant\max\{|t|,\tau\}\bigr\}$, and
	\begin{eqnarray}
		U_{t,\tau}(s)
		&\defeq& \id_t(s) -\id_0(s) -\tfrac{t}{2\tau} \bigl(\id_\tau(s)-\id_{-\tau}(s) \bigr),\quad s\in\R,\label{lm_UVW_def_U}\\
		V_{t,\tau}(s)&\defeq&\int_{-\infty}^s \rho_{t,\tau}(y) \dd y,\qquad   \rho_{t,\tau}(s)
		\defeq-\tfrac{1}{2}\Bigl(|s-t|-|s|-\tfrac{t}{2\tau}\,\bigl(|s-\tau|-|s+\tau|\bigr)\Bigr),\quad s\in\R,\label{lm_UVW_def_V_rho}\\
		W_{t,\tau}(s)&\defeq&U_{t,\tau}(s)+V_{t,\tau}(s),\quad s\in\R.\label{lm_UVW_def_W}
	\end{eqnarray}
	For any $t\in\R$ and $\tau>0$ it is true that $U_{t,\tau}(s)=0$ and $\rho_{t,\tau}(s)= 0$ for all $s\notin A_{t,\tau}$, $\rho_{t,\tau}$  is a continuous function on $\R$ with a broken-line graph, and, in particular, $\rho_{t,\tau}\in L_1(\R)$, the functions $U_{t,\tau}$, $V_{t,\tau}$, and $W_{t,\tau}$ belong to the class $\BV$.
\end{Lemma}
\textbf{Proof of Lemma \ref{lm_UVW}.} Let us fix $t\in\R$, $\tau>0$, and define $A_{t,\tau}$ as in the formulation. We write
\begin{eqnarray}\label{conc_U}
	e^{itx} -1 -\tfrac{it}{\tau} \sin(\tau x)=	e^{itx} -1 -\tfrac{t}{2\tau}\bigl(e^{i\tau x}-e^{-i\tau x}\bigr)=\int_{\R} e^{isx} \dd U_{t,\tau}(s),\quad x\in\R,
\end{eqnarray}
where $U_{t,\tau}$ is defined by \eqref{lm_UVW_def_U}. Using the definition of the function $\id_a(\cdot)$, $a\in\R$, it is easily seen that $U_{t,\tau}$ is an right-continuous function on $\R$, $U_{t,\tau}(s)=0$ for all $s\notin A_{t,\tau}$, and, in particular, $U_{t,\tau}\in \BV$.  Therefore the set $\R$  can be changed by $A_{t,\tau}$ in the integral \eqref{conc_U}. 

Let us consider the function
\begin{eqnarray*}
	\varphi_{t,\tau}(x)\defeq \Bigl(e^{itx} -1 -\tfrac{it}{\tau} \sin(\tau x)\Bigr)\dfrac{1}{x^2},\quad x\in\R.
\end{eqnarray*}
Observe that $\varphi_{t,\tau}\in L_1(\R)$. So we define
\begin{eqnarray}\label{def_rho_proof}
	\rho_{t,\tau}(s)\defeq\dfrac{1}{2\pi} \int_{\R} e^{-isx}\varphi_{t,\tau}(x)  \dd x,\quad s\in\R.
\end{eqnarray}
Let us find an explicit formula for $\rho_{t,\tau}(s)$ for every $s\in\R$. Observe that $x\mapsto \Realpart \varphi_{t,\tau}(x)$, $x\in\R$, is an even function and $x\mapsto \Imagpart \varphi_{t,\tau}(x)$, $x\in\R$, is an odd function. Therefore 
\begin{eqnarray*}
	\rho_{t,\tau}(s)&=&\dfrac{1}{2\pi} \int_{\R} \bigl(\Realpart \varphi_{t,\tau}(x)\, \cos (sx)+ \Imagpart \varphi_{t,\tau}(x)\, \sin (sx)\bigr)\dd x\\
	&=&\dfrac{1}{\pi} \int_{0}^\infty \bigl(\Realpart \varphi_{t,\tau}(x)\, \cos (sx)+ \Imagpart \varphi_{t,\tau}(x)\, \sin (sx)\bigr)\dd x\\
	&=& \dfrac{1}{\pi} \int_{0}^\infty \biggl(\dfrac{\cos(tx) -1 }{x^2}\, \cos(sx) + \dfrac{\sin(tx) -\tfrac{t}{\tau} \sin(\tau x)}{x^2}\, \sin(sx) \biggr)\dd x,\quad s\in\R.
\end{eqnarray*}
Next, using the known trigonometric formulas, we  write
\begin{eqnarray*}
	\rho_{t,\tau}(s)
	&=&\dfrac{1}{\pi} \int_{0}^\infty \biggl(\dfrac{\cos(tx)\cos(sx)+\sin(tx)\sin(sx) -\cos(sx) }{x^2}  - \dfrac{\tfrac{t}{\tau} \sin(\tau x)\sin(sx)}{x^2}\,  \biggr)\dd x\\
	&=&\dfrac{1}{\pi} \int_{0}^\infty \biggl(\dfrac{\cos((s-t)x) -\cos(sx) }{x^2}  - \dfrac{t}{\tau}\,\dfrac{ \cos((s-\tau)x)-\cos((s+\tau)x) }{2x^2}\,  \biggr)\dd x\\
	&=&\dfrac{1}{\pi} \int_{0}^\infty \dfrac{\cos(|s-t|x) -\cos(|s|x) }{x^2}\dd x  - \dfrac{t}{2\tau}\cdot\dfrac{1}{\pi} \int_{0}^\infty\dfrac{ \cos(|s-\tau|x)-\cos(|s+\tau|x) }{x^2}\dd x, \quad s\in\R.
\end{eqnarray*}
It is known (see \cite{GradRyz} p. 450, formula \textbf{3.782} 2.) that 
\begin{eqnarray*}
	\int_0^\infty \dfrac{1-\cos(ax)}{x^2}\, \dd x= \dfrac{a\pi}{2},\quad a\geqslant 0. 
\end{eqnarray*}
Hence
\begin{eqnarray*}
	\rho_{t,\tau}(s)
	=-\tfrac{1}{2}\Bigl(|s-t|-|s|-\tfrac{t}{2\tau}\,\bigl(|s-\tau|-|s+\tau|\bigr)\Bigr),\quad s\in\R,
\end{eqnarray*}
as in \eqref{lm_UVW_def_V_rho}. We see that $\rho_{t,\tau}$ is a continuous function with a broken-line graph. Also observe that $\rho_{t,\tau}(s)=0$ for all $s\notin A_{t,\tau}$. Indeed, if $s>\max\{|t|,\tau\}$, then
\begin{eqnarray*}
	\rho_{t,\tau}(s)=-\tfrac{1}{2}\Bigl(s-t-s-\tfrac{t}{2\tau} \bigl(s-\tau-(s+\tau)\bigr)\Bigr)=-\tfrac{1}{2}\Bigl(-t-\tfrac{t}{2\tau}\cdot (-2\tau)\Bigr)=-\tfrac{1}{2}(-t+t)=0,
\end{eqnarray*}
and if $s<-\max\{|t|,\tau\}$, then
\begin{eqnarray*}
	\rho_{t,\tau}(s)=-\tfrac{1}{2}\Bigl(-(s-t)+s-\tfrac{t}{2\tau} \bigl(-(s-\tau)+s+\tau\bigr)\Bigr)=-\tfrac{1}{2}\Bigl(t-\tfrac{t}{2\tau}\cdot 2\tau\Bigr)=-\tfrac{1}{2}(t-t)=0.
\end{eqnarray*}
Thus $\rho_{t,\tau}\in L_1(\R)$. By the way, observe that $V_{t,\tau}$, which is defined by \eqref{lm_UVW_def_V_rho}, is a  continuous function on $\R$ and it vanishes at $-\infty$, i.e. $V_{t,\tau}\in\BV$. Then, according to these remarks and \eqref{def_rho_proof}, we have
\begin{eqnarray*}
	\varphi_{t,\tau}(x)=\int_{\R} e^{isx}\rho_{t,\tau}(s)  \dd s=\int_{A_{t,\tau}} e^{isx}\rho_{t,\tau}(s)  \dd s=\int_{A_{t,\tau}} e^{isx}  \dd V_{t,\tau}(s),\quad x\in\R.
\end{eqnarray*}

Next, summing the proved equalities in \eqref{lm_UVW_intUV}, we get \eqref{lm_UVW_intW} with $W_{t,\tau}$ defined by \eqref{lm_UVW_def_W}. Since $U_{t,\tau}$ and $V_{t,\tau}$ belong to $\BV$, we conclude that $W_{t,\tau}\in\BV$.\quad $\Box$\\

\textbf{Proof of Theorem \ref{th_FnconvtoF}.}  Let $f$ be a characteristic function of the limit distribution function $F$. By the continuity theorem, we have 
\begin{eqnarray}\label{conc_fn_conv}
	f_n(t)\to f(t),\quad n\to\infty,\quad \text{for every}\,\,\, t\in\R.
\end{eqnarray}
Moreover, it is well known (see \cite{Lukacs}) that
\begin{eqnarray}\label{conc_fn_uniconv}
	\sup_{t\in[-T,T]}|f_n(t)-f(t)|\to 0,\quad n\to\infty,\quad \text{for any}\quad T>0.
\end{eqnarray}

First let us recall that characteristic functions of quasi-infinitely divisible distributions have no zeroes on the real line (see \cite{LindPanSato} or  \eqref{ineq_absfn0} below).  So, in particular, $f_n(t)\ne 0$, $t\in\R$, $n\in\N$. We now show that $f(t)\ne 0$ for all $t\in\R$.  For any fixed  $n\in\N$ and $t\in\R$ we consider
\begin{eqnarray*}
	|f_n(t)|&=& \exp\biggl\{\int_{\R} \bigl(\cos(tx) -1 \bigr)\tfrac{1+x^2}{x^2} \dd G_n(x)\biggr\}\\
	&\geqslant& \exp\biggl\{-\biggl|\int_{\R} \bigl(\cos(tx) -1 \bigr)\tfrac{1+x^2}{x^2} \dd G_n(x)\biggr|\biggr\}\\
	&\geqslant& \exp\biggl\{-\int_{\R} \bigl(1-\cos(tx) \bigr)\tfrac{1+x^2}{x^2} \dd |G_n|(x)\biggr\}.
\end{eqnarray*}
Let us estimate the inner function $x\mapsto \bigl(1-\cos(tx) \bigr)\tfrac{1+x^2}{x^2}$, $x\in\R$, which is equal to $\tfrac{t^2}{2}$ at $x=0$ for the continuity by the well known convention. Due to the inequality $1-\cos y\leqslant \tfrac{y^2}{2}$, $y\in\R$, for the case $|tx|\leqslant 2$ we have
\begin{eqnarray*}
	\bigl(1-\cos(tx) \bigr)\tfrac{1+x^2}{x^2}\leqslant \tfrac{t^2 x^2}{2}\cdot \tfrac{1+x^2}{x^2}=\tfrac{t^2+ t^2x^2}{2}=\tfrac{t^2}{2}+2.
\end{eqnarray*}
Using the simple inequalty $1-\cos y\leqslant  2$, $y\in\R$, for the case $|tx|>2$ we obtain
\begin{eqnarray*}
	\bigl(1-\cos(tx) \bigr)\tfrac{1+x^2}{x^2}\leqslant 2\cdot \tfrac{1+x^2}{x^2}= 2\cdot\bigl(\tfrac{1}{x^2}+1\bigr)   =2\cdot\bigl(\tfrac{t^2}{4}+1\bigr)=\tfrac{t^2}{2}+2.
\end{eqnarray*}
Thus 
\begin{eqnarray}\label{ineq_cos}
	\bigl(1-\cos(tx) \bigr)\tfrac{1+x^2}{x^2}\leqslant \tfrac{t^2}{2}+2, \quad\text{for any}\quad x\in\R,\quad t\in\R.  
\end{eqnarray}
Thus for any  $n\in\N$ and $t\in\R$ we obtain 
\begin{eqnarray}\label{ineq_absfn0}
	|f_n(t)|\geqslant \exp\biggl\{-\int_{\R} \Bigl(\tfrac{t^2}{2}+2\Bigr) \dd |G_n|(x)\biggr\}=\exp\Bigl\{-\Bigl(\tfrac{t^2}{2}+2\Bigr)   \|G_n\|\Bigr\}>0.
\end{eqnarray}
Hence, due to  \eqref{assum_UNB_Gn} and \eqref{conc_fn_conv}, we have
\begin{eqnarray*}
 |f(t)|=\lim_{n\to\infty} |f_n(t)|\geqslant \exp\Bigl\{-\Bigl(\tfrac{t^2}{2}+2\Bigr)\varlimsup_{n\to\infty}  \|G_n\|\Bigr\}=  \exp\Bigl\{-\Bigl(\tfrac{t^2}{2}+2\Bigr) B\Bigr\}>0,\quad t\in\R
\end{eqnarray*}
i.e. $f(t)\ne 0$ for any $t\in\R$.

Due to the above remarks,  the distinguished logarithms $t\mapsto\Ln f(t)$ and $t\mapsto \Ln f_n(t)$, $n\in\N$,  are defined for all $t\in\R$. According to \eqref{def_fn}, we have 
\begin{eqnarray}\label{def_Lnfn}
	\Ln f_n(t)= it \gamma_n+\int_{\R} \Bigl(e^{itx} -1 -\tfrac{it}{\tau} \sin(\tau x)\Bigr)\tfrac{1+x^2 }{x^2}  \dd G_n(x),\quad t\in\R,\quad n\in\N.
\end{eqnarray}
Due to the convergence \eqref{conc_fn_conv}, we have that
\begin{eqnarray}\label{conc_conv_Lnfnf}
	\Ln f_n(t)\to \Ln f(t),\quad n\to\infty,\quad \text{for every}\quad t\in\R.
\end{eqnarray}
Hence, in particular, 
\begin{eqnarray*}
	\gamma_n=\dfrac{\Imagpart (\Ln f_n(\tau))}{\tau}\to \dfrac{\Imagpart (\Ln f(\tau))}{\tau}\in\R,\quad n\to\infty.
\end{eqnarray*}
We denote this limit by $\gamma$. So we have
\begin{eqnarray}\label{conc_convgamman}
	\gamma_n\to \gamma,\quad n\to\infty. 
\end{eqnarray} 

We next introduce the following functions
\begin{eqnarray*}
	\psi(t,s)\defeq\Ln f(t)-\tfrac{1}{2} \bigl(\Ln f(t-s)+\Ln f(t+s)\bigr),\quad t\in\R,\quad s\geqslant 0,
\end{eqnarray*}
and analogously
\begin{eqnarray}\label{def_psin}
	\psi_n(t,s)\defeq\Ln f_n(t)-\tfrac{1}{2} \bigl(\Ln f_n(t-s)+\Ln f_n(t+s)\bigr),\quad t\in\R,\quad s\geqslant 0, \quad n\in\N.	 
\end{eqnarray}
From \eqref{conc_conv_Lnfnf} we conclude that
\begin{eqnarray}\label{conc_conv_psin}
	\psi_n(t,s)\to\psi(t,s),\quad n\to\infty, \quad \text{for any}\quad t\in\R,\quad s\geqslant 0. 
\end{eqnarray}
Moreover, since \eqref{conc_fn_uniconv} implies the convergence (see \cite{Khinch}, p. 15)
\begin{eqnarray*}
	\sup_{t\in[-T,T]}|\Ln f_n(t)-\Ln f(t)|\to 0,\quad n\to\infty,\quad \text{for any}\quad T>0,
\end{eqnarray*}
it is clear that 
\begin{eqnarray}\label{conc_psin_uniconv}
	\sup_{t,s\in[-T,T]}\bigl|\psi_n(t,s)-\psi(t,s)  \bigr|\to 0,\quad n\to\infty,\quad \text{for any}\quad T>0.
\end{eqnarray}

We next show that $\psi_n$, $n\in\N$, are uniformly bounded over $t\in\R$ and $n\in\N$ for any fixed $s\geqslant 0$. Using \eqref{def_Lnfn} in \eqref{def_psin}, we have
\begin{eqnarray}
	\psi_n(t,s)&=& it \gamma_n+\int_{\R} \Bigl(e^{itx} -1 -\tfrac{it}{\tau} \sin(\tau x)\Bigr)\tfrac{1+x^2 }{x^2}  \dd G_n(x)\nonumber\\
	&&{}- \tfrac{1}{2}\biggl(i2t \gamma_n+\int_{\R} \Bigl(e^{itx}\bigl(e^{-isx}+e^{isx}\bigr) -2 -\tfrac{i2t}{\tau} \sin(\tau x)\Bigr)\tfrac{1+x^2 }{x^2}  \dd G_n(x) \biggr)\nonumber\\
    &=&\int_{\R} e^{itx} \bigl(1-\cos(sx) \bigr) \tfrac{1+x^2}{x^2}\dd G_n(x),\quad t\in\R,\quad s\geqslant 0,\quad n\in\N.\label{repr_psin}
\end{eqnarray}
The estimate \eqref{ineq_cos} yields 
\begin{eqnarray*}
	\sup_{t\in\R}|\psi_n(t,s)|&\leqslant& \sup_{t\in\R}\int_{\R}  \Bigl|e^{itx}\bigl(1-\cos(sx) \bigr) \tfrac{1+x^2}{x^2}\Bigr|\dd |G_n|(x)\\
	&=&\int_{\R}  \bigl(1-\cos(sx) \bigr) \tfrac{1+x^2}{x^2}\dd |G_n|(x)\leqslant \Bigl(\tfrac{s^2}{2}+1\Bigr) \|G_n\|,\quad s\geqslant 0,\quad n\in\N.
\end{eqnarray*}
According to \eqref{assum_UNB_Gn}, there exists a constant $B_0\geqslant0$ such that  $\|G_n\|\leqslant B_0$ for all $n\in\N$. Then we conclude
\begin{eqnarray}\label{ineq_psin_s}
	\sup_{n\in\N}\sup_{t\in\R}|\psi_n(t,s)|\leqslant B_0\cdot\Bigl(\tfrac{s^2}{2}+1\Bigr),\quad s\geqslant0. 
\end{eqnarray}
Additionally, in view of \eqref{conc_conv_psin}, we obtain
\begin{eqnarray}\label{ineq_psi_s}
	\sup_{t\in\R}|\psi(t,s)|\leqslant B_0\cdot\Bigl(\tfrac{s^2}{2}+1\Bigr),\quad s\geqslant0. 
\end{eqnarray}

Next, since $\int_0^\infty (s^2+1) e^{-s} \dd s<\infty$, we can define the functions
\begin{eqnarray*}
	g_n(t)\defeq\int\limits_0^\infty \psi_n(t,s) e^{-s}\dd s,\quad n\in \N,\quad\text{and}\quad g(t)\defeq\int\limits_0^\infty  \psi(t,s) e^{-s}\dd s,\quad t\in\R,
\end{eqnarray*}
and, due to \eqref{conc_conv_psin}, conclude at once by the Lebesgue dominated convergence theorem that
\begin{eqnarray*}
	g_n(t)\to g(t),\quad n\to\infty, \quad \text{for every}\,\,\, t\in\R.
\end{eqnarray*}
Let us prove that
\begin{eqnarray}\label{conc_gn_uniconv}
	\sup_{t\in[-T,T]} |g_n(t)-g(t)|\to 0,\quad n\to\infty, \quad\text{for any}\quad T>0.
\end{eqnarray}
We fix any $T>0$ and $\e>0$. It is clear that for every $n\in\N$
\begin{eqnarray}\label{ineq_gn_psin}
	\sup_{t\in[-T,T]} |g_n(t)-g(t)|\leqslant  \int\limits_0^\infty \sup_{t\in[-T,T]}\bigl|\psi_n(t,s)-\psi(t,s)\bigr|\, e^{-s}\dd s.
\end{eqnarray}
We denote by  $J_n(T)$ the last integral for every $n\in\N$. Let us choose a constant $h_\e>0$ such that 
\begin{eqnarray}\label{ineq_int_s2exp_e}
	B_0\int\limits_{h_\e}^\infty (s^2+2) e^{-s} \dd s <\e.
\end{eqnarray}
Then we write $J_n(T)=J_{n,1}(T)+J_{n,2}(T)$, $n\in\N$, where
\begin{eqnarray*}
	J_{n,1}(T)\defeq   \int\limits_0^{h_\e} \sup_{t\in[-T,T]}\bigl|\psi_n(t,s)-\psi(t,s)\bigr|\, e^{-s}\dd s,\qquad J_{n,2}(T)\defeq   \int\limits_{h_\e}^{\infty} \sup_{t\in[-T,T]}\bigl|\psi_n(t,s)-\psi(t,s)\bigr|\, e^{-s}\dd s.	
\end{eqnarray*}
All the integrals $J_n(T)$, $J_{n,1}(T)$, and $J_{n,2}(T)$ are non-negative. Observe that
\begin{eqnarray*}
	J_{n,1}(T)\leqslant   \sup_{\substack{t\in[-T,T],\\ s\in[0,h_\e]}}\bigl|\psi_n(t,s)-\psi(t,s)\bigr|\int\limits_0^{h_\e}  e^{-s}\dd s\leqslant \sup_{t,s\in[-T_\e,T_\e]}\bigl|\psi_n(t,s)-\psi(t,s)\bigr|,\quad n\in\N,
\end{eqnarray*}
where $T_\e\defeq \max\{T, h_\e\}$. Due to \eqref{conc_psin_uniconv},  the last supremum vanishes as $n\to\infty$. So there exists $n_\e\in\N$ such that $J_{n,1}(T)<\e$ for any $n\geqslant n_\e$. Let us turn to $J_{n,2}(T)$. According to \eqref{ineq_psin_s}, \eqref{ineq_psi_s}, and \eqref{ineq_int_s2exp_e}, we have
\begin{eqnarray*}
	J_{n,2}(T)\leqslant \int\limits_{h_\e}^{\infty}\Bigl(\, \sup_{t\in\R}\bigl|\psi_n(t,s)\bigr|+\sup_{t\in\R}\bigl|\psi(t,s)\bigr|\Bigr)\, e^{-s}\dd s\leqslant \int\limits_{h_\e}^\infty B_0\bigl(s^2+2\bigr)\, e^{-s} \dd s <\e.
\end{eqnarray*}
Then $J_n(T)=J_{n,1}(T)+J_{n,2}(T)<2\e$ for any $n\geqslant n_\e$. Since $\e>0$ was chosen arbitrarily, $J_n(T)\to 0$ as $n\to\infty$. Thus, according to \eqref{ineq_gn_psin}, we obtain \eqref{conc_gn_uniconv}. Since $g$ is a uniform limit  of continuous functions $g_n$ on any segment $[-T,T]$ as $n\to\infty$, the function $g$ is continuous on $\R$.

Let us consider the functions $g_n$, $n\in\N$. Using \eqref{repr_psin}, we write
\begin{eqnarray*}
	g_n(t)&=& \int_0^\infty \psi_n(t,s) e^{-s}\dd s\\
	&=&\int_0^\infty\biggl(\,\int_{\R} e^{itx} \bigl(1-\cos(sx) \bigr) \tfrac{1+x^2}{x^2}\dd G_n(x)\biggr)e^{-s}\dd s\\
	&=&\int_{\R}   \biggl(\,\int_{0}^{\infty}\bigl(1-\cos(sx) \bigr)e^{-s}\dd s\biggr)\,e^{itx}\,\tfrac{1+x^2}{x^2} \dd G_n(x),\quad t\in\R,\quad n\in\N.
\end{eqnarray*}
The inner integral is calculated (see \cite{GradRyz} p. 486, formula \textbf{3.893} 2.):
\begin{eqnarray*}
	\int_{0}^{\infty}\bigl(1-\cos(sx) \bigr)e^{-s}\dd s=1- \int_{0}^{\infty}\cos(sx) e^{-s}\dd s=1-\dfrac{1}{1+x^2}=\dfrac{x^2}{1+x^2},\quad x\in\R.
\end{eqnarray*}
Therefore we have
\begin{eqnarray*}
	g_n(t)=\int_{\R} e^{itx} \dd G_n(x),\quad t\in\R, \quad n\in\N.
\end{eqnarray*}

We now use Theorem \ref{th_FStoFBV}. So there exists a unique function $G\in\BV$  such that $\|G\|\leqslant B$ and the equality
\begin{eqnarray*}
	g(t)=\int_{\R} e^{itx} \dd G(x)
\end{eqnarray*}
holds \textit{for all} $t\in\R$, because $g$ is continuous on $\R$. Moreover, due to the theorem, we have $G_n\cconv G$ and also $G_n(+\infty)\to G(+\infty)$, $n\to\infty$. 

We now prove that for any $t\in\R$ and $\tau>0$ 
\begin{eqnarray}\label{conc_conv_intGnG}
	\int_{\R} \Bigl(e^{itx} -1 -\tfrac{it}{\tau} \sin(\tau x)\Bigr)\tfrac{1+x^2 }{x^2}  \dd G_n(x)\to \int_{\R} \Bigl(e^{itx} -1 -\tfrac{it}{\tau} \sin(\tau x)\Bigr)\tfrac{1+x^2 }{x^2}  \dd G(x),\quad n\to\infty.
\end{eqnarray}
Let us fix $t\in\R$ and $\tau>0$. From Lemma \ref{lm_UVW} we know that
\begin{eqnarray*}
	\Bigl(e^{itx} -1 -\tfrac{it}{\tau} \sin(\tau x)\Bigr)\tfrac{1+x^2 }{x^2}= \int_{A_{t,\tau}} e^{isx} \dd W_{t,\tau}(s),\quad x\in\R,
\end{eqnarray*}
where $A_{t,\tau}=\bigl\{s\in\R: |s|\leqslant \max\{|t|,\tau\}\bigr\}$, and $W_{t,\tau}\in\BV$. Hence for every $n\in\N$
\begin{eqnarray*}
	\int_{\R} \Bigl(e^{itx} -1 -\tfrac{it}{\tau} \sin(\tau x)\Bigr)\tfrac{1+x^2 }{x^2}  \dd G_n(x)&=&\int_{\R} \biggl(\int_{A_{t,\tau}} e^{isx} \dd W_{t,\tau}(s)\biggr) \dd G_n(x)\\
	&=& \int_{A_{t,\tau}}\biggl(\int_{\R} e^{isx}\dd G_n(x)\biggr)  \dd W_{t,\tau}(s)\\
	&=&\int_{A_{t,\tau}}g_n(s)\dd W_{t,\tau}(s).
\end{eqnarray*}
Also we have analogously that
\begin{eqnarray*}
	\int_{\R} \Bigl(e^{itx} -1 -\tfrac{it}{\tau} \sin(\tau x)\Bigr)\tfrac{1+x^2 }{x^2}  \dd G(x)=\int_{A_{t,\tau}}g(s)\dd W_{t,\tau}(s).
\end{eqnarray*}
Thus \eqref{conc_conv_intGnG} takes the form
\begin{eqnarray*}
	\int_{A_{t,\tau}}g_n(s)\dd W_{t,\tau}(s)\to\int_{A_{t,\tau}}g(s)\dd W_{t,\tau}(s),\quad n\to\infty.
\end{eqnarray*}
This convergence holds. Indeed, for every $n\in\N$
\begin{eqnarray*}
	\biggl|\int_{A_{t,\tau}}g_n(s)\dd W_{t,\tau}(s)-\int_{A_{t,\tau}}g(s)\dd W_{t,\tau}(s)\biggr|&\leqslant& \int_{A_{t,\tau}}\bigl|g_n(s)\dd W_{t,\tau}(s)-g(s)\bigr|\dd |W_{t,\tau}|(s)\\
	&\leqslant& \sup_{s\in A_{t,\tau}}\bigl|g_n(s)-g(s)\bigr|\cdot \|W_{t,\tau}\|,
\end{eqnarray*}
where, due to \eqref{conc_gn_uniconv}, the supremum vanishes as $n\to\infty$. Thus we proved \eqref{conc_conv_intGnG}.

From \eqref{def_Lnfn}, \eqref{conc_convgamman}, and  \eqref{conc_conv_intGnG}, for any $t\in\R$ we have
\begin{eqnarray*}
	\Ln f_n(t)&=& it \gamma_n+\int_{\R} \Bigl(e^{itx} -1 -\tfrac{it}{\tau} \sin(\tau x)\Bigr)\tfrac{1+x^2 }{x^2}  \dd G_n(x)\\
	&\to& it \gamma+\int_{\R} \Bigl(e^{itx} -1 -\tfrac{it}{\tau} \sin(\tau x)\Bigr)\tfrac{1+x^2 }{x^2}  \dd G(x),\quad n\to\infty.
\end{eqnarray*}
According to \eqref{conc_conv_Lnfnf}, we conclude that
\begin{eqnarray*}
	\Ln f(t)=it \gamma+\int_{\R} \Bigl(e^{itx} -1 -\tfrac{it}{\tau} \sin(\tau x)\Bigr)\tfrac{1+x^2 }{x^2}  \dd G(x),\quad t\in\R,
\end{eqnarray*}
where, as we have already proved, $\gamma\in\R$ and $G\in\BV$. Thus $f$ has the L\'evy-Khinchine type representation with $(\gamma, G)$, i.e. the distribution function $F$ corresponding to  $f$ is quasi-infinitely divisible. \quad $\Box$\\

\textbf{Proof of Theorem \ref{th_FnconvtoFminus}.}  Let $f$ be a characteristic function of the limit distribution function $F$. So we have \eqref{conc_fn_conv} and also \eqref{conc_fn_uniconv} (see comments in the proof of Theorem \ref{th_FnconvtoF}).

Recall that $f_n(t)\ne 0$, $t\in\R$, $n\in\N$. Let us choose $\delta>0$ such that $f(t)\ne 0$, $|t|\leqslant \delta$ (it is possible because $f$ is continuous on $\R$ and $f(0)=1$). Let us consider values of the Khinchine functional $\chi_\delta(\cdot)$ (see \cite{LinOstr} p. 79) with parameter $\delta$ on $f$ and $f_n$, $n\in\N$:
\begin{eqnarray*}
	\chi_\delta(f)=- \dfrac{1}{\delta}\int_{0}^{\delta} \ln |f(s)|\dd s,\qquad \chi_\delta(f_n)=- \dfrac{1}{\delta}\int_{0}^{\delta} \ln |f_n(s)|\dd s,\quad n\in\N.
\end{eqnarray*}
These quatities are finite and nonegative. Due to \eqref{conc_fn_uniconv}, we have 
\begin{eqnarray}\label{conc_chiffn}
	\chi_\delta(f_n)\to \chi_\delta(f),\quad n\to\infty.
\end{eqnarray}
Observe that
\begin{eqnarray}
	\chi_\delta(f_n)&=&-\dfrac{1}{\delta} \int_0^\delta \biggl( \int_{\R}\bigl(\cos(sx)-1\bigr) \tfrac{1+x^2}{x^2}\dd G_n(x) \biggr) \dd s\nonumber\\
	&=&\int_{\R}\biggl(\dfrac{1}{\delta} \int_0^\delta  \bigl(1-\cos(sx)\bigr) \dd s\biggr)\tfrac{1+x^2}{x^2}\dd G_n(x)\nonumber\\
	&=&\int_{\R}\Bigl(1-\tfrac{\sin (\delta x)}{\delta x}\Bigr)\tfrac{1+x^2}{x^2}\dd G_n(x),\quad n\in\N.\label{conc_chideltfn}
\end{eqnarray}
where we set
\begin{eqnarray}\label{conv_cos_sin}
	\bigl(\cos(sx)-1\bigr) \tfrac{1+x^2}{x^2}\Bigr|_{x=0}=-\tfrac{s^2}{2},\qquad  \Bigl(1-\tfrac{\sin (\delta x)}{\delta x}\Bigr)\tfrac{1+x^2}{x^2}\Bigr|_{x=0}=\tfrac{\delta^2}{3!},
\end{eqnarray}
according to known expansions $\cos y=1-\tfrac{y^2}{2}+o(y^2)$ and $\sin y=y-\tfrac{y^3}{3!}+o(y^3)$, $y\to 0$. Let us consider the inner function of the integral in \eqref{conc_chideltfn}:
\begin{eqnarray*}
	x\mapsto \Bigl(1-\tfrac{\sin (\delta x)}{\delta x}\Bigr)\tfrac{1+x^2}{x^2},\quad x\in\R.
\end{eqnarray*}
By convention \eqref{conv_cos_sin}, it is continuous at the point $x=0$. We see that  this function is continuous and strictly positive on $\R$. Also observe that it tends to $1$ as $x\to\pm \infty$. Hence it is clear that there exist positive constants $c_\delta$ and $C_\delta$ such that
\begin{eqnarray}\label{ineq_sin}
	0<c_\delta\leqslant\Bigl(1-\tfrac{\sin (\delta x)}{\delta x}\Bigr)\tfrac{1+x^2}{x^2}\leqslant C_\delta<\infty,\quad x\in\R.
\end{eqnarray}

Let us take some decompositions \eqref{eq_Gn_decom} for $G_n$, $n\in\N$. According to \eqref{conc_chideltfn}, we have 
\begin{eqnarray*}
	\chi_\delta(f_n)=\int_{\R}\Bigl(1-\tfrac{\sin (\delta x)}{\delta x}\Bigr)\tfrac{1+x^2}{x^2}\dd G_n^+(x)-\int_{\R}\Bigl(1-\tfrac{\sin (\delta x)}{\delta x}\Bigr)\tfrac{1+x^2}{x^2}\dd G_n^-(x),\quad n\in\N.
\end{eqnarray*}
Due to \eqref{ineq_sin}, we obtain
\begin{eqnarray*}
	\chi_\delta(f_n)\geqslant c_\delta\int_{\R}\dd G_n^+(x)-C_\delta\int_{\R}\dd G_n^-(x)=c_\delta G_n^+(+\infty)-C_\delta G_n^-(+\infty),\quad n\in\N.
\end{eqnarray*}
From this we have
\begin{eqnarray*}
	 G_n^+(+\infty)\leqslant \dfrac{\chi_\delta(f_n)+C_\delta G_n^-(+\infty)}{c_\delta},\quad n\in\N.
\end{eqnarray*}
Hence, due to \eqref{assum_GnNeg} and \eqref{conc_chiffn}, we get
\begin{eqnarray}\label{conc_Gnplus}
	\varlimsup_{n\to\infty} G_n^+(+\infty)\leqslant \dfrac{1}{c_\delta}\bigl( \varlimsup_{n\to\infty} \chi_\delta(f_n)+ C_\delta\varlimsup_{n\to\infty}G_n^-(+\infty)\bigr)=\dfrac{1}{c_\delta}\bigl(  \chi_\delta(f)+ C_\delta M\bigr)<\infty.
\end{eqnarray}
According to \eqref{eq_Gn_decom} and conventions there, it is true that
\begin{eqnarray}\label{ineq_varGn}
	\|G_n\|\leqslant \|G_n^+\|+\|G_n^-\|= G_n^+(+\infty)+G_n^-(+\infty),\quad n\in\N.
\end{eqnarray}
So we conclude from  \eqref{assum_GnNeg} and \eqref{conc_Gnplus}  that \eqref{assum_UNB_Gn} holds for some $B<\infty$.

Thus all assertions of Theorem \ref{th_FnconvtoF} hold. In particular, $G_n\cconv G$ and $G_n(+\infty)\to G(+\infty)$, $n\to\infty$, where $G$ is some function from $\BV$.  It remains to prove that $B\leqslant G(+\infty)+2M$. Using inequality \eqref{ineq_varGn}, we write
\begin{eqnarray*}
	B=\varlimsup_{n\to\infty}\|G_n\|\leqslant \varlimsup_{n\to\infty}\bigl(G_n^+(+\infty)+G_n^-(+\infty)\bigr) \leqslant\varlimsup_{n\to\infty} \bigl(G_n^+(+\infty)-G_n^-(+\infty) \bigr)+2 \varlimsup_{n\to\infty} G_n^-(+\infty),
\end{eqnarray*}
but $ G_n(+\infty)=G_n^+(+\infty)-G_n^-(+\infty)$, $n\in\N$, and we obtain
\begin{eqnarray*}
	B\leqslant\varlimsup_{n\to\infty} G_n(+\infty)+2 \varlimsup_{n\to\infty} G_n^-(+\infty)=G(+\infty)+2 M,
\end{eqnarray*}
as required.\quad $\Box$\\

\textbf{Proof of Theorem \ref{th_FnconvtoF_suff}.}  Let $(F_{n_k})_{k\in\N}$ be arbitrary subsequence of $(F_{n})_{n\in\N}$, which weakly converges to some distribution function $F_*$. Due to the assumption of relative compactness of $(F_{n})_{n\in\N}$, such subsequence exists. By Theorem \ref{th_FnconvtoF}, $F_*$ is quasi-infinitely divisible with some spectral pair $(\gamma_*, G_*)$, where $\gamma_*\in\R$ and $G_*\in\BV$.   Moreover,  $\gamma_{n_k}\to \gamma_*$ and $G_{n_k}\cconv G_*$, $k\to\infty$. According to the assumption that $\gamma_n\to \gamma$, $n\to\infty$, we conclude that $\gamma_*=\gamma$. Let us show that $G_*=G$. By definition, the convergence  $G_{n_k}\cconv G_*$, $k\to\infty$, implies an existence of a subsequence $(G_{n'_l})_{l\in\N}$ in $(G_{n_k})_{k\in\N}$ such that
\begin{eqnarray*}
	G_{n'_l}(x_2)-G_{n'_l}(x_1) \to G_*(x_2)-G_*(x_1),\quad l\to\infty,
\end{eqnarray*}
for all $x_1,x_2\in\R$ except at most countable set $D'$. Due to the assumption that $G_n\cconv G$, $n\to\infty$, we can choose a further subsequence $(G_{n''_l})_{l\in\N}$ in $(G_{n'_l})_{l\in\N}$ such that
\begin{eqnarray*}
	G_{n''_l}(x_2)-G_{n''_l}(x_1) \to G(x_2)-G(x_1),\quad l\to\infty,
\end{eqnarray*}
for all $x_1,x_2\in\R$ except at most countable set $D''$ (and let $D'\ne D''$ in general).  Therefore
\begin{eqnarray*}
	G_*(x_2)-G_*(x_1) = G(x_2)-G(x_1)
\end{eqnarray*}
for all $x_1,x_2\in\R$ except at most countable set $D'\cup D''$. Letting $x_1\to -\infty$ over $x_1\in\R\setminus\!(D'\cup D'')$ we have $G_*(x_2)=G(x_2)$ for every $x_2\in \R\setminus\!(D'\cup D'')$ and, consequently, for all $x_2\in\R$, because $G_*, G \in \BV$, i.e. they are right-continuous and  $G_*(-\infty)=G(-\infty)=0$. Thus we proved that $\gamma_*=\gamma$ and $G_*=G$.  

The previous remark means that $(\gamma, G)$ is the spectral pair for some quasi-infinitely divisible distribution function $F$. We also saw that every  subsequence $(F_{n_k})_{k\in\N}$, which weakly converges to some distribution function, converges exactly to $F$, because a spectral pair uniquely determines a distribution function. Therefore, since $(F_n)_{n\in\N}$ is relatively compact, we conclude that whole  sequence $(F_n)_{n\in\N}$ weakly converges to $F$ (this is known fact, see  \cite{Bill} p. 337).   \quad $\Box$

\section{Acknowledgments}

This research was supported by the Ministry of Science and Higher Education of the Russian Federation, agreement 075-15-2019-1620 date 08/11/2019 and 075-15-2022-289 date 06/04/2022.

\end{document}